# Modified empirical CLT's under only pre-Gaussian conditions

Shahar Mendelson[1],[*] and Joel Zinn[2],[†]

*ANU & Technion I.I.T and Texas A&M University*

**Abstract:** We show that a modified Empirical process converges to the limiting Gaussian process whenever the limit is continuous. The modification depends on the properties of the limit via Talagrand's characterization of the continuity of Gaussian processes.

## 1. Introduction

Given a class of functions $F \subset L_2(\mu)$, a now standard method to use (iid) data to uniformly estimate or predict the mean of one of the functions in the class, is through the use of empirical processes. One has to bound the random variable

$$\sup_{f \in F} \left| \frac{1}{n} \sum_{i=1}^n f(X_i) - \mathbb{E}f \right| \equiv \|P_n - P\|_F.$$

One possibility that comes to mind is to use the fact that there is a bounded Gaussian process indexed by $F$ to bound $\sqrt{n}\|P_n - P\|_F$. To illustrate the difficulty one encounters, note that to use the classical Central Limit Theorem in finite dimensions, one only needs finite variance or, equivalently, the existence of the (Gaussian) limit. However, when working in the infinite dimensional situation there are sometimes non-trivial side conditions other than just the set being pregaussian. Those are connected to the random geometry of the set $F$ that are needed to ensure the existence of a useful bound. For example, one such situation is when the class is a collection of indicators of sets. If such a class does not satisfy the VC (Vapnik-Červonenkis) condition, then, in addition to knowing that the limiting Gaussian process is continuous, one has to check, for example, the following:

$$\frac{\log \#\{C \cap \{X_1, \ldots, X_n\} : C \in \mathcal{C}\}}{\sqrt{n}} \to 0 \text{ in outer probability} .$$

In this note we try to get around this problem by looking at a variant of the standard empirical process for which only the existence of the limiting Gaussian process is required to obtain both tail estimates and the Central Limit Theorem for the modified process.

The motivation for our study were the articles [7, 8], which focus on the following problem. Consider a density $p(x)$ on $\mathbb{R}$ which has a support that is a finite union of intervals, and on this support $p$ satisfies that $c^{-1} \leq p(x) \leq c$ and $|p(x) - p(y)| \leq c|x - y|$. It was shown in [7, 8] that under this assumption there is a histogram rule $\tilde{P}_n$ for which the following holds. If $F$ is a $P$-pregaussian class of indicator functions, or if $F$ is a $P$-pregaussian class of functions bounded by 1 which satisfy a certain metric entropy condition, then $\sqrt{n}(\tilde{P}_n - P)$ converges weakly to the limiting Gaussian process.

[1]The Australian National University, Canberra, ACT 0200, Australia and, Technion - Israel Institute of Technology, Haifa 32000, Israel, e-mail: shahar.mendelson@anu.edu.au
[2]Department of Mathematics, Texas A&M University, College Station, TX 77843-3386, e-mail: jzinn@math.tamu.edu
[*]Research partially supported by an Australian Research Council Discovery Grant DP0343616.
[†]Research partially supported by NSA Grant H98230-04-1-0108.

It seems that the method used in [7, 8] can not be extended to densities in $\mathbb{R}^2$, and even in the one-dimensional case the convergence result holds for a very restricted set of densities. Thus, our aim is to find some other empirical estimator for which the existence of the limiting Gaussian would imply convergence as above. Our estimator is based on Theorem 1 in [9] (see also Theorem 1.3 [4] and Theorem 14.8 [6]).

We begin with several facts and notation we will use throughout this note. If $G$ is a (centered) Gaussian process indexed by the set $T$, then for every $s, t \in T$, $\rho_2(s, t) = \left(\mathbb{E}(G_s - G_t)^2\right)^{1/2}$.

To prove that a process (and in particular, the modified empirical process we define) converges to the Gaussian process indexed by $F$, we require the notion of asymptotic equicontinuity.

**Definition 1.1.** A net $X_\alpha : \Omega_\alpha \longrightarrow \ell_\infty(T)$ is asymptotically uniformly $\rho$-equicontinuous in probability, if for every $\epsilon, \eta > 0$, there exists a $\delta > 0$ such that

$$\limsup_\alpha \Pr^*(\sup_{\rho(s,t)<\delta} |X_\alpha(s) - X_\alpha(t)| > \epsilon) < \eta.$$

**Theorem 1.2** ([13], 1.5.10). *Let $G$ be a Gaussian process and let $X_\alpha$ be a net of random variables with values in $\ell_\infty(T)$. Then there exists a version of $G$ which is a tight Borel measurable map into $\ell_\infty(T)$, and $X_\alpha$ converges weakly to $G$ if and only if*

  (i) *the finite dimensional distributions of $X_\alpha$ converge weakly to the corresponding finite dimensional distributions of $G$,*
 (ii) *$X_\alpha$ is asymptotically uniformly equicontinuous in probability with respect to $\rho_2$,*
(iii) *$(T, \rho_2)$ is totally bounded.*

The main technical tool we require is *generic chaining* which was developed by Talagrand (see [11] for the most recent survey on this topic).

**Definition 1.3.** For a metric space $(T, d)$, an *admissible sequence* of $T$ is a collection of subsets of $T$, $\{T_s : s \geq 0\}$, such that for every $s \geq 1$, $|T_s| = 2^{2^s}$ and $|T_0| = 1$. For $p \geq 1$, we define the $\gamma_p$ functional by

$$\gamma_p(T, d) = \inf \sup_{t \in T} \sum_{s=0}^\infty 2^{s/p} d(t, T_s),$$

where the infimum is taken with respect to all admissible sequences of $T$.

For every integer $s$ define the function $\pi_s : T \to T_s$, which maps every $t \in T$ to a nearest element to $t$ in $T_s$.

Using the $\gamma_p$ functionals one can bound the supremum of a process which satisfies an increment condition. As an example, recall the well known Bernstein inequality.

**Theorem 1.4.** *There exists an absolute constant c for which the following holds. Let $(\Omega, \mathcal{S}, P)$ be a probability space, let $f \in L_\infty(P)$ and let $X_1, \ldots, X_n$ be independent random variables distributed according to $P$. Then*

$$Pr\left(\left\{\left|\frac{1}{n}\sum_{i=1}^n f(X_i) - \mathbb{E}f\right| \geq t\right\}\right) \leq 2\exp\left(-cn\min\left\{\frac{t^2}{\|f\|_{L_2}^2}, \frac{t}{\|f\|_{L_\infty}}\right\}\right).$$

When one combines Bernstein's inequality with the generic chaining method, the following corollary is evident:

**Corollary 1.5** ([11])**.** *There exists an absolute constant c for which the following holds. Let $F$ be a class of functions on a probability space $(\Omega, P)$. Then, for every integer $n$,*

$$\mathbb{E}\|P_n - P\|_F \leq c\left(\frac{\gamma_2(F, \|\cdot\|_2)}{\sqrt{n}} + \frac{\gamma_1(F, \|\cdot\|_\infty)}{n}\right).$$

A key result, which follows from the generic chaining method, is that the expectation of the supremum of the Gaussian process indexed by $F$ and has a covariance structure endowed by $L_2(P)$ is finite if and only if $\gamma_2(F, L_2(P))$ is finite. Moreover, the result is quantitative in nature.

**Theorem 1.6.** *There exist absolute constants $c_1$ and $c_2$ for which the following holds. Let $F$ be a class of functions on $(\Omega, P)$, and let $G$ be the Gaussian process indexed by $F$. Then,*

$$c_1\gamma_2(F, L_2(P)) \leq \mathbb{E}\sup_{f \in F}|G_f| \leq c_2\gamma_2(F, L_2(P)).$$

The upper bound is due to Fernique [3], while the lower bound in due to Talagrand [10]. A proof of both parts can be found in [11]. From here on we denote $\mathbb{E}\sup_f |G_f|$ by $\mathbb{E}\|G\|_F$.

In a similar fashion one can formulate a continuity condition for the Gaussian process indexed by $F$ using the generic chaining machinery.

**Theorem 1.7** ([11] Theorem 1.4.1)**.** *Consider a Gaussian process, $\{G_f : f \in F\}$, where $F$ is countable. Then the following are equivalent:*

(1) *The map $f \longrightarrow G_f(\omega)$ is uniformly continuous on $(F, \|\cdot\|_2)$ with probability one.*
(2) *We have*
$$\lim_{\epsilon \to 0} \mathbb{E}\sup_{\|G_f - G_{f'}\|_2 \leq \epsilon}|G_f - G_{f'}| = 0.$$
(3) *There exist an admissible sequence of partitions of $F$ such that*
$$\lim_{s_0 \to \infty} \sup_{f \in F} \sum_{s \geq s_0} 2^{s/2}d(f, F_s) = 0.$$

Note that the admissible sequence in (3) can be taken as an almost optimal admissible sequence in the definition of $\gamma_2(F, \|\cdot\|_2)$, at a price of an absolute constant. Indeed, by combining the two admissible sequences $(T_s)$ and $(T'_s)$ - the first - an almost optimal one from the definition of $\gamma_2$ and the second from (1.7), we obtain a new admissible sequence for which $T_s \cup T'_s \subset T''_{s+1}$. Thus, we may assume that the almost optimal sequence in (3) satisfies that $\sup_{f \in F} \sum_{s=0}^\infty 2^{s/2}d_{\|\cdot\|_2}(t, T_s) \sim \mathbb{E}\|G\|_F$.

Finally, a notational convention. Throughout, all absolute constants will be denoted by $c$, $C$ or $K$. Their values may change from line to line. We write $a \sim b$ if there are absolute constants $c$ and $C$ such that $ca \leq b \leq Ca$.

## 2. The main theorems

Let $F \subset L_2(P)$ and assume that there exists a Gaussian process indexed by $F$ such that $\mathbb{E}\|G\|_F < \infty$. By Theorem 1.6, $\gamma_2(F, L_2(P)) \sim \mathbb{E}\|G\|_F$, and let $(F_s)_{s=0}^\infty$, $F_s \subset F$ be an almost optimal admissible sequence with respect to the $L_2(P)$ norm as described above, that is,

$$\sup_{f \in F} \sum_{s=0}^{\infty} 2^{s/2} \|f - \pi_s(f)\|_{L_2(P)} \leq c\gamma_2(F, L_2(P)),$$

and

$$\lim_{s_0 \to \infty} \sup_{f \in F} \sum_{s \geq s_0} 2^{s/2} \|f - \pi_s(f)\|_{L_2(P)} = 0.$$

Set $\Delta_s(f) = \pi_s(f) - \pi_{s-1}(f)$. For every $s$, $f$ and $\lambda \geq 0$ consider a truncated part of $\Delta_s(f)$, defined by $\Delta'_s(f, \lambda) = \Delta_s(f)\mathbb{1}_{\{|\Delta_s(f)| \leq \lambda\}}$. As will be made clear shortly, the truncation level $\lambda$ depends both on the specific $f \in F$, as well as on the size of the sample $n$ and the index $s$.

**Lemma 2.1.** *There exist absolute constants $c_1$, $c_2$ and $c_3$ for which the following holds. Let $F$ be a class of functions on the probability space $(\Omega, P)$, and set $X_1, \ldots, X_n$ to be independent, distributed according to $P$. Let $(F_s)_{s=0}^\infty$ be an admissible sequence in $F$ with respect to $L_2(P)$, and for every $f \in F$ set $\lambda = c_1\sqrt{n}\|\Delta_s(f)\|_2/2^{s/2}$. Then, for every $u > 1/2$ and every integer $s$, with probability at least $1 - 2\exp(-c_2 2^s \min\{u^2, u\})$, for every $f \in F$.*

$$(2.1) \qquad \left| \frac{1}{n} \sum_{i=1}^{n} (\Delta'_s(f, \lambda))(X_i) - \mathbb{E}\Delta'_s(f, \lambda) \right| \leq c_3 u \frac{2^{s/2}\|\Delta_s(f)\|_2}{\sqrt{n}}.$$

*Proof.* Let $c_1$ and $c_2$ be constants to be named later. By Bernstein's inequality [2], for $t = c_1 2^{s/2} u \|\Delta_s(f)\|_2/\sqrt{n}$, it is evident that for any $\lambda > 0$,

$$Pr\left(\left|\frac{1}{n}\sum_{i=1}^{n}(\Delta'_s(f,\lambda))(X_i) - \mathbb{E}\Delta'_s(f,\lambda)\right| \geq t\right)$$
$$\leq 2\exp\left(-cn\min\left\{\frac{t^2}{\|\Delta'_s(f,\lambda)\|_2^2}, \frac{t}{\lambda}\right\}\right).$$

Since $\|\Delta'_s(f,\lambda)\|_2 \leq \|\Delta_s(f)\|_2$ and $\lambda = c_2\sqrt{n}\|\Delta_s(f)\|_2/2^{s/2}$, then for the choice of $t$, it follows that with probability at least $1 - 2\exp(-c_3 2^s \min\{u^2, u\})$,

$$\left|\frac{1}{n}\sum_{i=1}^{n}(\Delta'_s(f,\lambda))(X_i) - \mathbb{E}\Delta'_s(f,\lambda)\right| \leq c_4 u \frac{2^{s/2}\|\Delta_s(f)\|_2}{\sqrt{n}},$$

where $c_3$ and $c_4$ depend only on $c_1$ and $c_2$. Thus, for an appropriate choice of these constants and since $|\{\Delta_s(f) : f \in F\}| \leq |F_s| \cdot |F_{s-1}| \leq 2^{2^{s+1}}$ the claim follows. □

Note that there is nothing magical with the lower bound of $1/2$ on $u$. Any other absolute constant would do, and would lead to changed absolute constants $c_1$, $c_2$ and $c_3$.

Using the Lemma we can define a process $\Phi_n$ for which $\|P_n(\Phi_n) - P\|_F \leq c\mathbb{E}\|G\|_F/\sqrt{n}$, and thus, the fact that the limiting Gaussian process exists suffices to yield a useful bound on the way in which the empirical estimator approximates the mean.

**Definition 2.2.** Let $\lambda(f,n,s) \equiv \lambda = c_0\sqrt{n}\|\Delta_s(f)\|_2/2^{s/2}$, where $c_0$ was determined in Lemma 2.1, and for each $s_0$ let

$$\Phi_{n,s_0}(f) = \sum_{s=s_0+1}^{\infty} \Delta'_s(f,\lambda)$$

and

$$\Phi_n(f) = \sum_{s=1}^{\infty} \Delta'_s(f,\lambda)$$

**Theorem 2.3.** *There exist absolute constants $c_1$ and $c_2$ for which the following holds. Let $F$ and $X_1, \ldots, X_n$ be as above. Then, the mapping $\Phi_n : F \to L_1(P)$, and for every $u > 1/2$, with probability at least*

$$1 - 2\exp(-c_1\min\{u^2, u\}),$$

*for every $f \in F$*

$$(2.2) \qquad \left|\frac{1}{n}\sum_{i=1}^{n}(\Phi_n(f))(X_i) - \mathbb{E}f\right| \leq c_2(u+1)\frac{\mathbb{E}\|G\|_F}{\sqrt{n}}.$$

*and also with probability at least $1 - 2\exp(-c_1 2^{s_0}\min\{u^2, u\})$*

$$(2.3) \qquad \sup_{f \in F}\left|\frac{1}{n}\sum_{i=1}^{n}\Big(\Phi_{n,s_0}(f)(X_i) - \mathbb{E}\Phi_{n,s_0}(f)\Big)\right| \leq \frac{c_2 u}{\sqrt{n}}\sup_{f \in F}\sum_{s=s_0+1}^{\infty} 2^{s/2}\|\Delta_s(f)\|_2.$$

*Proof.* Without loss of generality, assume that $0 \in F$ and that $\pi_0 = 0$. Let $(F_s)_{s=1}^{\infty}$ be an almost optimal admissible sequence, and in particular, by Theorem 1.6

$$\sup_{f \in F}\sum_{s=1}^{\infty} 2^{s/2}\|\Delta_s(f)\|_2 \leq K\mathbb{E}\|G\|_F$$

for a suitable absolute constant $K$.

Note that as $\pi_0(f) = 0$ for every $f \in F$ one can write $f = \sum_{s=1}^{\infty}\Delta_s(f)$. Let us show that $\Phi_n$, and therefore $\Phi_{n,s_0}$, are well defined and maps $F$ into $L_1(P)$. Indeed, since $\sum_{s=1}^{\infty}\Delta_s(f)$ converges in $L_2(P)$, it suffices to prove that $\sum_{s=1}^{\infty}(\Delta_s(f) - \Delta'_s(f,\lambda)) \equiv \sum_{s=1}^{\infty}\Delta''_s(f)$ converges in $L_1(P)$. Observe that for every $f \in F$,

$$\mathbb{E}|\Delta''_s(f,\lambda)| = \mathbb{E}|\Delta_s(f)|\mathbb{1}_{\{|\Delta_s(f)|>\lambda\}}$$

$$(2.4) \qquad \leq \|\Delta_s(f)\|_2\left(Pr\left(|\Delta_s(f)| > \lambda\right)\right)^{1/2} \leq \frac{\|\Delta_s(f)\|_2^2}{\lambda}$$

$$\leq \frac{2^{s/2}\|\Delta_s(f)\|_2}{c_0\sqrt{n}}.$$

Since $\sum_{s=1}^{\infty} 2^{s/2}\|\Delta_s(f)\|_2$ converges for every $f$, it implies that $\Phi_n$ is well defined and takes values in $L_1$.

By Lemma 2.1, with probability at least

$$(2.5) \qquad 1 - 2\sum_{s=s_0+1}^{\infty}\exp(-c_1 2^s\min\{u^2, u\}) \geq 1 - 2\exp(-c_2 2^{s_0}\min\{u^2, u\}),$$

$$\sup_{f \in F} \sum_{s=s_0+1}^{\infty} \left| \frac{1}{n} \sum_{i=1}^{n} \left(\Delta'_s(f,\lambda)\right)(X_i) - \mathbb{E}\Delta'_s(f,\lambda) \right|$$

$$\leq \frac{c_3 u}{\sqrt{n}} \sup_{f \in F} \sum_{s=s_0+1}^{\infty} 2^{s/2} \|\Delta_s(f)\|_2,$$

and when $s_0 = 0$ we'll use that by Theorem 1.6 this last quantity is

$$\leq c_4 u \frac{\mathbb{E}\|G\|_F}{\sqrt{n}}.$$

Hence, with that probability, for every $f \in F$,

$$\left| \frac{1}{n} \sum_{i=1}^{n} (\Phi_n(f)(X_i) - \mathbb{E}f) \right| \leq \left| \frac{1}{n} \sum_{i=1}^{n} (\Phi_n(f)(X_i) - \mathbb{E}\Phi_n(f)) \right| + |\mathbb{E}f - \mathbb{E}\Phi_n(f)|$$

$$\leq c_5 u \frac{\mathbb{E}\|G\|_F}{\sqrt{n}} + \left| \mathbb{E} \sum_{s=1}^{\infty} \Delta''_s(f) \right| \leq c_6(u+1) \frac{\mathbb{E}\|G\|_F}{\sqrt{n}},$$

where the last term is estimated using the same argument as in (2.4) and the inequality (2.2) in Theorem 2.3.

We also have that with the promised lower bound on the probability,

$$\sup_{f \in F} \left| \frac{1}{n} \sum_{i=1}^{n} \left(\Phi_{n,s_0}(f)(X_i) - \mathbb{E}\Phi_{n,s_0}(f)\right) \right| \leq \frac{c_3 u}{\sqrt{n}} \sup_{f \in F} \sum_{s=s_0+1}^{\infty} 2^{s/2} \|\Delta_s(f)\|_2. \qquad \square$$

Next, we prove a limit theorem for $\{\sqrt{n}(P_n - P)(\Phi_n)(f) : f \in F\}$ and show that we can replace $\mathbb{E}\Phi_n(f)$ with $\mathbb{E}f$ and still obtain a limit theorem. For this we need to prove an inequality for the oscillation of the process $\sqrt{n}(P_n - P)(\Phi_n(f))$. To that end, define $Q_n := \sqrt{n}(P_n - P)$.

**Proposition 2.4.** *Let $F$ be a class of functions on $(\Omega, P)$, such that the Gaussian process indexed by $F$ exists and is continuous. If $\Phi_n$ is as above, then for any $\eta > 0$,*

$$\lim_{\delta \to 0} \lim_{n \to \infty} Pr\left(\sup_{\|f-\tilde{f}\|_2 < \delta} \left|Q_n(\Phi_n(f) - \Phi_n(\tilde{f}))\right| > \eta\right) = 0.$$

*Proof.* By the definition of $\Phi_n$ which uses an almost optimal admissible sequence, for every $\delta > 0$ there is some $s_0$ such that

$$\sup_{f \in F} \sum_{s=s_0}^{\infty} 2^{s/2} \|f - \pi_s(f)\|_2 \leq \delta,$$

hence for any $f, \tilde{f} \in F$, $\|\pi_{s_0}(f) - \pi_{s_0}(\tilde{f})\|_2 < 2\delta + \|f - \tilde{f}\|_2$. Using the notation of Theorem 2.3, put $\Phi_{n,s_0}(f) := \sum_{s=s_0+1}^{\infty} \Delta'_s(f, \lambda)$ and $\Psi_{n,s_0} = \Phi_n - \Phi_{n,s_0}$. Therefore,

$$\begin{aligned} I : &= Pr\left(\sup_{\|f-\tilde{f}\|_2<\delta} \left|Q_n\left(\Phi_n(f) - \Phi_n(\tilde{f})\right)\right| > \eta\right) \\ &= Pr\Big(\sup_{\|f-\tilde{f}\|_2<\delta} \Big|Q_n\left(\Psi_{n,s_0}(f) - \Psi_{n,s_0}(\tilde{f})\right) \\ &\quad + Q_n(\Phi_{n,s_0}(f)) - Q_n(\Phi_{n,s_0}(\tilde{f}))\Big| > \eta\Big) \\ &\leq Pr\left(\sup_{\|\pi_{s_0} f - \pi_{s_0}\tilde{f}\|_2<3\delta} \left|Q_n\left(\Psi_{n,s_0}(f) - \Psi_{n,s_0}(\tilde{f})\right)\right| > \frac{\eta}{3}\right) \\ &\quad + 2Pr\left(\sup_f |Q_n(\Phi_{n,s_0}(f))| > \frac{\eta}{3}\right) \\ &:= (II) + (III) \end{aligned}$$

From the proof of Theorem 2.3 by integrating tail probabilities

$$\mathbb{E}\sup_{f\in F} |Q_n\left(\Phi_{n,s_0}(f)\right)| \leq c \sup_{f\in F} \sum_{s>s_0} 2^{s/2} \|f - \pi_s(f)\|_2$$

which by Theorem 1.7(3) and our choice of the admissible sequence converges to 0 as $s_0 \to \infty$. Furthermore, by the finite dimensional Central Limit Theorem $\lim_{\delta\to 0} \lim_{n\to\infty} (II) = 0$, which completes the proof. □

Hence, we know that $Q_n$ is asymptotically uniformly equicontinuous. We'll now prove the other necessary ingredients needed to show that $Q_n$ converges to the original Gaussian process.

**Proposition 2.5.** *Let $F$ and $\Phi_n$ be as in Proposition 2.4. Then the following holds:*

(i) $\lim_{n\to\infty} \sqrt{n} \sup_{f\in F} |\mathbb{E}\Phi_n(f) - \mathbb{E}f| = 0$,
(ii) *For every $f \in F$*, $\lim_{n\to\infty} \|\Phi_n(f) - f\|_2 = 0$,
(iii) *For every $f \in F$*, $\lim_{n\to\infty} \mathbb{E}\max_{j\leq n} \dfrac{|\Phi_n(f)(X_j)|^2}{n} = 0$.

*Proof.*

1. Let $s_0$ be an integer to be named later and set

$$\lambda = c_0 \sqrt{n} \|\Delta_s(f)\|_2 / 2^{s/2}$$

as in Lemma 2.1. In particular, the set $\{\Delta_s(f) : s \leq s_0, f \in F\}$ is finite and for every $f \in F$,

$$\sqrt{n}\mathbb{E}|\Delta_s(f)|\mathbb{1}_{\{|\Delta_s(f)|>\lambda\}} \leq \frac{2^{s/2}}{c_0\|\Delta_s(f)\|_2} \mathbb{E}|\Delta_s(f)|^2 \mathbb{1}_{\{|\Delta_s(f)|>\lambda\}}$$

which, by the definition of $\lambda$ tends to 0 as $n$ tends to infinity. Hence, for every fixed $s_0$,

$$\lim_{n\to\infty} \sum_{s=1}^{s_0} \sqrt{n}\mathbb{E}|\Delta_s(f)|\mathbb{1}_{\{|\Delta_s(f)|>\lambda\}} = 0.$$

Therefore, for every $s_0$,

$$\lim_{n\to\infty} \sqrt{n} \sup_{f\in F} |\mathbb{E}\Phi_n(f) - \mathbb{E}f|$$
$$\leq \lim_{n\to\infty} \sqrt{n} \sup_{f\in F} \left( \sum_{s=1}^{s_0} \mathbb{E}|\Delta_s(f)|\mathbb{1}_{\{|\Delta_s(f)|>\lambda\}} + \sum_{s>s_0} \mathbb{E}|\Delta_s(f)|\mathbb{1}_{\{|\Delta_s(f)|>\lambda\}} \right)$$
$$\leq c_2 \sup_{f\in F} \sum_{s>s_0} 2^{s/2} \|\Delta_s(f)\|_2,$$

where the last inequality is evident from (2.4) and the choice of a suitable absolute constant $c_2$.

2. Again, we shall use the fact that for every fixed $f$ and $s$, $\lambda$ depends on $n$ and tends to 0 as $n$ tends to infinity. Clearly, for every fixed $s_0$,

$$\|f - \Phi_n(f)\|_2 \leq \sum_{s\leq s_0} \|\Delta_s(f)\mathbb{1}_{\{|\Delta_s(f)|>\lambda\}}\|_2 + \sum_{s>s_0} \|\Delta_s(f)\|_2$$
$$\leq \sum_{s\leq s_0} \|\Delta_s(f)\mathbb{1}_{\{|\Delta_s(f)|>\lambda\}}\|_2 + c_3\gamma_2(F, L_2) \sum_{s>s_0} 2^{-s/2}$$

For an absolute constant $c_3$. Indeed, this follows from the fact that for every $s$,

$$2^{s/2}\|\Delta_s(f)\|_2 \leq \sum_{s=0}^{\infty} 2^{s/2}\|\Delta_s(f)\|_2 \leq c_3\gamma_2(F, L_2),$$

and of course the constant $c_3$ does not depend on $s$.

Hence, for every fixed $f \in F$,

$$\limsup_{n\to\infty} \|f - \Phi_n(f)\|_2 \leq c_3\gamma_2(F, L_2) \sum_{s>s_0} 2^{-s/2}$$

for every $s_0$, and this last quantity goes to zero as $s_0 \to \infty$.

3. If $f(X)$ is square integrable then for any $b > 0$,

$$\limsup_{n\to\infty} \mathbb{E}\max_{j\leq n} \frac{|f(X_j)|^2}{n} \leq \limsup_{n\to\infty} \mathbb{E}\left(\frac{b^2}{n} + \frac{1}{n}\sum_{j\leq n} f^2(X_j)\mathbb{1}_{\{|f(X_j)|>b\}}\right)$$
$$= \mathbb{E}f^2(X)\mathbb{1}_{\{|f(X)|>b\}}.$$

Since the left hand side does not depend on $b$ and the right hand side converges to zero as $b$ tends to $\infty$, $\lim_{n\to\infty} \mathbb{E}\max_{j\leq n} \frac{|f(X_j)|^2}{n} = 0$. Therefore, to complete the proof it suffices to show that

$$\mathbb{E}\max_{j\leq n} \frac{|f(X_j) - \Phi_n(f)(X_j)|^2}{n} \to 0.$$

But, using (2),

$$\mathbb{E}\max_{j\leq n} \frac{|f(X_j) - \Phi_n(f)(X_j)|^2}{n} \leq \frac{1}{n}\sum_{j\leq n} \mathbb{E}|f(X_j) - \Phi_n(f)(X_j)|^2$$
$$= \mathbb{E}|(f - \Phi_n(f))(X)|^2 \to 0. \qquad \square$$

The final ingredient we require is the following result on triangular arrays.

**Lemma 2.6.** *For each $n$, let $\{\xi_{n,j}\}_{j=1}^{n}$ by nonnegative, square integrable, independent random variables for which $\lim_{n\to\infty} \mathbb{E}\max_{j\leq n} \xi_{n,j}^2 = 0$. Then, for every $\delta > 0$, $\lim_{n\to\infty} \sum_{j=1}^{n} \mathbb{E}\xi_{n,j}^2 \mathbb{1}_{\{\xi_{n,j}\geq\delta\}} = 0$.*

*Proof.* Consider the stopping times

$$\tau = \tau_n = \begin{cases} \inf\{k \leq n : \xi_{n,k} > \delta\} & \text{if } \max_{r\leq n} \xi_{n,r} > \delta \\ \infty & \text{if } \max_{r\leq n} \xi_{n,r} \leq \delta. \end{cases}$$

Then, (see [12]) for every $n$

$$\mathbb{E}\max_{j\leq n}\xi_{n,j}^2 \geq \mathbb{E}\xi_{n,\tau_n}^2 \mathbb{1}_{\{\tau_n<\infty\}} = \sum_{l=1}^{n} \mathbb{E}\xi_{n,l}^2 \mathbb{1}_{\{\xi_{n,l}>\delta,\max_{i\leq l-1}\xi_{n,i}\leq\delta\}}$$

$$= \sum_{l=1}^{n} \mathbb{E}\xi_{n,l}^2 \mathbb{1}_{\{\xi_{n,l}>\delta\}} \Pr(\max_{i\leq l-1} \xi_{n,i} \leq \delta)$$

$$\geq \sum_{l=1}^{n} \mathbb{E}\xi_{n,l}^2 \mathbb{1}_{\{\xi_{n,l}>\delta\}} \Pr(\max_{i\leq n} \xi_{n,i} \leq \delta)$$

The result now follows, since the hypothesis implies that this last probability converges to one as $n$ tends to infinity. □

We now can conclude

**Theorem 2.7.** *If the Gaussian process indexed by $F$ is continuous then $\{\sqrt{n}(P_n(\Phi_n(f)) - Pf) : f \in F\}$ converges to $\{G_f : f \in F\}$.*

*Proof.* By Theorem 1.2 and Proposition 2.4 we only need to show that

(i) the finite dimensional distributions of $Q_n$ converge to those of $G$ and
(ii) $(F, \rho_2)$ is totally bounded.

For (i) we need to check that for any $\{f_i\}_{i=1}^{k} \subseteq F$,

$$(Q_n(\Phi_n(f_1)), \ldots, Q_n(\Phi_n(f_k)))$$

converges in distribution to $(G_{f_1}, \ldots, G_{f_k})$. To see this we apply the Cramer-Wold device, that is, by noting that to show the convergence in distribution, we only have to check that the characteristic function (on $\mathbb{R}^k$) converges, and hence it suffices to show that any finite linear combination of $\{Q_n(\Phi_n(f_i))\}_{i=1}^{k}$, say, $\sum_{i=1}^{k} a_i Q_n(\Phi_n(f_i))$ converges in distribution to $\sum_{i=1}^{k} a_i G_{f_i}$. To verify this, recall the classical Central Limit Theorem for triangular arrays (see, e.g, [5] or [1] Theorem 3.5). Namely, it suffices to prove that

(a) for any $\eta > 0$, $\lim_{n\to\infty} \Pr(\max_{j\leq n} |\sum_{i=1}^{k} a_i \Phi_n(f_i(X_j))| > \eta) = 0$ and
(b) $\lim_{n\to\infty} \text{Var}((\sum_{i=1}^{k} a_i \Phi_n(f_i))\mathbb{1}_{\{|\sum_{i=1}^{k} a_i \Phi_n(f_i)|>\eta\}}) = \text{Var}(\sum_{i=1}^{k} a_i f_i)$.

(a) follows from Proposition 2.5(iii) and (ii) and (b) follows from 2.5(ii) and Lemma 2.6.
(ii) follows from the assumed continuity of $\{G_f : f \in F\}$ with respect to $\rho_2$ (see p. 41 [13]). □

## 3. Changing the level of truncation

The question we wish to tackle here is whether it is possible to find different "universal" truncation levels instead of $\sqrt{n}$, and still have a process $\Psi$ which is tight, and satisfies that $n^{-1}\sum_{i=1}^{n}\left(\psi_n(f)\right)(X_i)$ uniformly approximates $\mathbb{E}f$ (that is, one can replace $\mathbb{E}\Psi_n(f)$ with $\mathbb{E}f$). We show that such a uniform level of truncation has to be asymptotically larger than $\sqrt{n}$.

**Definition 3.1.** Given a class of functions $F$ and a non-decreasing sequence of positive numbers, $b = \{b_n\}_{n=1}^{\infty}$, let

$$\Phi_{n,b} = \sum_{s=1}^{\infty} \Delta_s(f) \mathbb{1}_{\{|\Delta_s(f)| \leq b_n \|\Delta_s(f)\|_2 / 2^{s/2}\}}.$$

**Definition 3.2.** A sequence of processes $\{U_n(f) : f \in F\}$ is said to be stochastically bounded if for every $\epsilon > 0$ there is a constant $C < \infty$ such that

$$\Pr(\sup_{f \in \mathcal{F}} |U_n(f)| > C) < \epsilon.$$

**Theorem 3.3.** *Assume that $\{b_n\}_n$ is an increasing sequence of positive numbers and that the probability space, $(\Omega, \mathcal{S}, P)$, is continuous. Assume also that for every pregaussian class of functions on $(\Omega, \mathcal{S}, P)$, the process*

$$\{\sup_{f \in F} \sqrt{n} |P_n(\Phi_{n,b}(f)) - \mathbb{E}f|\}_n$$

*is stochastically bounded. Then, there exists $\delta > 0$ such that $\inf_n \dfrac{b_n}{\sqrt{n}} > \delta$.*

*Proof.* Clearly, if $\{\sqrt{n}\left(P'_n(\Phi_{n,b}(f)) - \mathbb{E}f\right) : f \in F\}$ is based on an independent copy, $\{X'_j\}$, then $\{\sup_{f \in F} \sqrt{n} |P'_n(\Phi_{n,b}(f)) - \mathbb{E}f|\}_{n=1}^{\infty}$ is also stochastically bounded. Hence, the difference is stochastically bounded, and thus,

$$\{\sup_{f \in F} \sqrt{n} |P_n(\Phi_{n,b}(f)) - \mathbb{E}\Phi_{n,b}(f)|\}_n$$

is stochastically bounded, implying that $\sqrt{n} \sup_{f \in F} |\mathbb{E}f - \mathbb{E}\Phi_{n,b}(f)|$ is bounded.

In particular, for every nonnegative $f \in L_2(P)$, if we let $F = \{f, 0\}$, then the sequence $\{\sqrt{n}|\mathbb{E}f - \mathbb{E}\Phi_{n,b}(f)|\}_{n=1}^{\infty}$ is bounded. Note that in this case we may assume that $\pi_s(f) = f$ for $s \geq 1$ and $\pi_0(f) = 0$, implying that $\sqrt{n}\mathbb{E}f\mathbb{1}_{\{f > b_n \|f\|_2/\sqrt{2}\}}$ is bounded.

Observe that this implies that $\sqrt{n}\mathbb{E}f\mathbb{1}_{\{f > b_n\}}$ is bounded. Indeed, choose $b_{k_0}$ such that $\|f\mathbb{1}_{\{f > b_{k_0}\}}\|_2 \leq \sqrt{2}$. Applying the above to the function $h = f\mathbb{1}_{\{f > b_{k_0}\}}$, it follows that $\|h\|_2 \leq \sqrt{2}$ and

$$\sqrt{n}\mathbb{E}h\mathbb{1}_{\{h > b_n\}} = \sqrt{n}\mathbb{E}f\mathbb{1}_{\{f > b_{k_0}\}}\mathbb{1}_{\{f\mathbb{1}_{\{f > b_{k_0}\}} > b_n\}} = \sqrt{n}\mathbb{E}f\mathbb{1}_{\{f > b_{\max(k_0,n)}\}}.$$

Hence, $\sqrt{n}\mathbb{E}f\mathbb{1}_{\{f > b_n\}}$ is bounded, as claimed.

For every sequence $\{a_k\}_k$ for which $\sum_k |a_k|/k < \infty$, consider a function $f$ with $\Pr(f = b_k) = b_1^2 \dfrac{|a_k|}{kb_k^2}$. Such a function exists by the continuity of the probability space $(\Omega, \mathcal{S}, P)$. Then,

$$\mathbb{E}f^2 = \sum_k b_k^2 b_1^2 \dfrac{|a_k|}{kb_k^2} < \infty.$$

Therefore, $\mathbb{E}f\mathbb{1}_{\{f>b_k\}} = \sum_{l>k} b_l b_1^2 \frac{|a_l|}{lb_l^2}$, implying that for every sequence $\{a_k\}_{k=1}^\infty$ as above, $\sup_k \sqrt{k} \sum_{l>k} \frac{|a_l|}{l} < \infty$.

Consider the Banach spaces $B_1$ and $B_2$, endowed with the norms $\|\{a_k\}\|_1 := \sum_{k=1}^\infty \frac{|a_k|}{k}$ and $\|\{a_k\}\|_2 := \sup_{k\geq 1} \sqrt{k} \sum_{l>k} \frac{|a_l|}{lb_l}$. Note that the identity map $\mathbb{I}: B_1 \longrightarrow B_2$ is bounded using the Closed Graph Theorem. Indeed, for $A_n := \{a_{n,k}\}_{k=1}^\infty, B := \{b_k\}_{k=1}^\infty$ and $C := \{c_k\}_{k=1}^\infty$ assume that $\|A_n - B\|_1 \to 0$ and $\|A_n - C\|_2 \to 0$. These conditions respectively imply convergence coordinate-wise, that is, for every $r$, $\lim_{n\to\infty} a_{n,r} = b_r$ and $\lim_{n\to\infty} a_{n,r} = c_r$. Thus, $B = C$, and the graph is closed, implying that the map is bounded.

Therefore, there exists a constant, $C$, such that

$$(3.1) \qquad \sup_{k\geq 1} \sqrt{k} \sum_{l>k} \frac{|a_l|}{lb_l} \leq C \sum_{k=1}^\infty \frac{|a_k|}{k}.$$

Applying (3.1) to the sequence for which the $n^{\text{th}}$ term is one and others zero shows that for $n > 1$:

$$\frac{\sqrt{n-1}}{nb_n} \leq C\frac{1}{n},$$

from which the claim follows. $\square$

## References

[1] ARAUJO, A. AND GINÉ, E. (1980). *The Central Limit Theorem for Real and Banach Valued Random Variables*. Wiley Series in Probability and Mathematical Statistics. John Wiley & Sons, New York-Chichester-Brisbane. MR576407 (83e:60003)

[2] BENNET, G. (1962). Probability inequalities for sums of independent random variables. *JASA* **57** 33–45.

[3] FERNIQUE, X. (1975). *Regularité des trajectoires des fonctions aléatoires gaussiennes*. École d'Été de Probabilités de Saint-Flour, IV-1974. *Lectures Notes in Math.*, Vol. 480. Springer, Berlin, pp. 1–96. MR0413238 (54 #1355)

[4] GINÉ, E. AND ZINN, J. (1986). Lectures on the central limit theorem for empirical processes. Probability and Banach spaces (Zaragoza, 1985). *Lectures Notes in Math.*, Vol. 1221. Springer, Berlin, pp. 50–113. MR88i:60063

[5] GNEDENKO, B. V. AND KOLMOGOROV, A. N. (1968). *Limit Distributions for Sums of Independent Random Variables*. Translated from the Russian, annotated, and revised by K. L. Chung. With appendices by J. L. Doob and P. L. Hsu. Revised edition, Addison-Wesley Publishing Col., Reading, MA–London–Don Mills., Ont. MR0233400 (38 #1722)

[6] LEDOUX, M. AND TALAGRAND, M. (1991). *Probability in Banach Spaces*. Ergebnisse der Mathematik und ihrer Grenzgebiete (3) [Results in Mathematics and Related Areas (3) Isoperimetry and processes], Vol. 23. Springer-Verlag, Berlin. MR1102015 (93c:60001)

[7] RADULOVIĆ, D. AND WEGKAMP, M. (2000). Weak convergence of smoothed empirical processes: beyond Donsker classes. *High Dimensional Probability*, II (Seattle, WA, 1999), *Progr. Probab.*, Vol. 47. Birkhäuser Boston, Boston, MA, pp. 89–105. MR1857317 (2002h:60043)

[8] RADULOVIĆ, D. AND WEGKAMP, M. (2003). Necessary and sufficient conditions for weak convergence of smoothed empirical processes. *Statist. Probab. Lett.* **61** (3) 321–336. MR2003i:60041

[9] TALAGRAND, M. (1987). Donsker classes and random geometry. *Ann. Probab.* **15** (4) 1327–1338. MR89b:60090
[10] TALAGRAND, M. (1987). Regularity of Gaussian processes. *Acta Math.* **159** (1–2) 99–149. MR906527 (89b:60106)
[11] TALAGRAND, M. (2005). *The Generic Chaining. Upper and Lower Bounds of Stochastic Processes.* Springer Monographs in Mathematics. Springer-Verlag, Berlin. MR2133757
[12] VAKHANIA, N. N., TARIELADZE, V. I. AND CHOBANYAN, S. A. (1987). *Probability Distributions on Banach spaces.* Mathematics and Its Applications (Soviet Series), Vol. 14. D. Reidel Publishing Co., Dordrecht. Translated from the Russian and with a preface by Wojbor A. Woyczynski. MR1435288 (97k:60007)
[13] VAN DER VAART, AAD W. AND WELLNER, JON A. (1996). *Weak Convergence and Empirical Processes. With Applications to Statistics.* Springer Series in Statistics. Springer-Verlag, New York. MR1385671 (97g:60035)